\documentclass{jhrs}

\usepackage{amsmath}

\input xy

\xyoption{all}

\def\diag{\ensuremath{\,\mathrm{diag}\,}}

\newtheorem{theorem}{Theorem}

\newtheorem{lemma}[theorem]{Lemma}

\begin{document}

\title{Diagonal fibrations are pointwise fibrations\\\vspace{0.3cm}
{\sl \small Dedicated To The Memory Of Saunders Mac Lane(1909-2005)}}
\thanks{Supported by DGI of Spain:
(MTM2004-01060 and MTM2006-06317), FEDER, and Junta de Andaluc\'ia: P06-FQM-1889.
}

\author{A.M. Cegarra}\email{acegarra@ugr.es}
\address{Departamento de \'{A}lgebra\\ Facultad de Ciencias\\ Universidad de
Granada\\ 18071 Granada\\ Spain}

\author{J. Remedios}\email{jremed@ull.es}
\address{Departamento de Matem\'{a}tica Fundamental\\ Universidad de La
Laguna\\ 38271 La Laguna\\ Spain}



\begin{abstract}
On the category of bisimplicial sets there are different Quillen closed model structures
associated to various definitions of fibrations. In one of them,  which is due to Bousfield and
Kan and that consists of seeing a bisimplicial set as a simplicial object in the category of
simplicial sets, fibrations are those bisimplicial set maps such that each of the induced
simplicial set maps is a Kan fibration, that is, the pointwise fibrations. In another of them,
introduced by Moerdijk, a bisimplicial map is a fibration if it induces a Kan fibration of
associated diagonal simplicial sets, that is, the diagonal fibrations. In this note, we prove that
every diagonal fibration is a pointwise fibration.
\end{abstract}

\received{Day Month Year}   
\revised{Day Month Year}    
\published{Day Month Year}  
\submitted{Name of Editor}  

\volumeyear{2007} 
\volumenumber{1}  
\issuenumber{2}   

\startpage{1}     

\maketitle
\section{Introduction and summary} There are several (Quillen) closed model structures on the
category of bisimplical sets, see \cite[IV, \S 3]{g-j}. This paper concerns two of them,
namely, the so-called Bousfield-Kan and Moerdijk structures, that we briefly recall below:

On the one hand, in the closed model structure by Bousfield-Kan, bisimplicial sets are regarded as
diagrams of simplicial sets and then fibrations are the pointwise Kan fibrations and weak
equivalences are the pointwise weak homotopy equivalences. To be more precise, a bisimplicial set
$X:\Delta^{op}\times \Delta^{op}\to \mbox{Set}$, $([p],[q])\mapsto X_{p,q}$, is seen as a
``horizontal" simplicial object in the category of ``vertical" simplicial sets, $X:\Delta^{op}\to
\mathbf{S}$, $[p]\mapsto X_{p,*}$ and then, a bisimplicial map $f:X\to Y$ is a fibration (resp. a
weak equivalence) if all simplicial maps $f_{p,*}:X_{p,*}\to Y_{p,*}$, $p\geq 0$, are Kan
fibrations (resp. weak homotopy equivalences).

On the other hand, the Moerdijk closed model structure on the bisimplicial set category is
transferred from the ordinary model structure on the simplicial set category  through the diagonal
functor, $X\mapsto \diag X: [n]\mapsto X_{n,n}$. Thus, in this closed model structure, a
bisimplicial map $f:X\to Y$ is a fibration (resp. a weak equivalence) if the induced diagonal
simplicial map $\diag f:\diag X\to \diag Y$ is a Kan fibration (resp. a weak homotopy
equivalence).

Several useful relationships between these two different homotopy theories of bisimplicial sets have
been established and, perhaps, the best known of them is the following:

\vspace{0.2cm} \noindent {\bf Theorem. }(Bousfield-Kan) {\em Let $f:X\to Y$ be a bisimplicial map
such that $f_{p,*}:X_{p,*}\to Y_{p,*}$ is a weak homotopy equivalence for each $p\geq 0$. Then
$\diag f:\diag X\to \diag Y$ is a weak homotopy equivalence.}

\vspace{0.2cm}The purpose of this brief note is to state and prove a suitable counterpart to
Bousfield-Kan's theorem for fibrations, namely:

\vspace{0.2cm}\noindent {\bf Theorem 1}. {\em Let $f:X\to Y$ be a bisimplicial map such that
$\diag f:\diag X\to \diag Y$ is a Kan fibration. Then $f_{p,*}:X_{p,*}\to Y_{p,*}$ is a Kan
fibration for each $p\geq 0$.}

\vspace{0.2cm}Note that the converse of Theorem 1 is not true in general. A counterexample is
given in the last section of the paper.

\vspace{0.3cm}{\bf Acknowledgements.} The authors are much indebted to the referees, whose useful observations greatly improved
our exposition. The second author is grateful to the Algebra
Department in the University of Granada for the excellent atmosphere and hospitality.

\section{Some preliminaries}
We use the standard conventions and terminology which can be found in texts on simplicial
homotopy theory, e. g. \cite{g-j} or \cite{may}. For definiteness or emphasis we state the
following.

We denote by $\Delta$  the category of finite ordered sets of integers $[n]=\{0,1,\ldots,n\}$,
$n\geq 0$, with weakly order-preserving maps between them. The category of {\em simplicial sets}
is the category of functors $X:\Delta^{op}\to \mbox{Set}$, where Set is the category of sets. If
$X$ is a simplicial set and $\alpha:[m]\to[n]$ is a map in $\Delta$, then we write
$X_n=X[n]$ and $\alpha^*=X(\alpha):X_n\to X_m$. Recall that all maps in $\Delta$ are
generated by the injections $\delta_i:[n-1]\to[n]$ (cofaces), $0\leq i\leq n$, which miss
out the $i$th element and the surjections $\sigma_i:[n+1]\to [n]$ (codegeneracies), $0\leq i\leq
n$, which repeat the $i$th element (see \cite[VII, \S5, Proposition 2]{mac}).  Thus, in order to
define a simplicial set, it suffices to give the  sets of {\em $n$-simplices} $X_n$, $n\geq 0$,
together with maps
$$
\begin{array}{lll}d_i=\delta_i^*:X_n\to X_{n-1},&\ 0\leq i\leq n& \ (\mbox{the face maps}),\\[3pt]
s_i=\sigma_i^*:X_n\to X_{n+1},&\ 0\leq i\leq n& \ (\mbox{the degeneracy maps}),
\end{array}
$$

satisfying the well-known basic {\em simplicial identities} such as $d_id_j=d_{j-1}d_i$ if $i<j$, etc.
(see \cite[p. 175]{mac}). In addition, we shall write down a list of other identities between
some iterated compositions of face and degeneracy maps, which will be used latter.
The proof of these equalities is straightforward and left to the reader.

\begin{lemma}\label{l1} On any simplicial set, the following equalities hold:
\begin{eqnarray}
\label{I1}d_id_j^m=&d_{j}^md_{i+m}\hspace{0.5cm}&\ \mbox{if\qquad} i\geq j\\[3pt]
\label{I2}d_i^m=&d_i^{m-1}d_j\hspace{0.5cm}&\ \mbox{if\qquad} i\leq j<i+m\\[3pt]
\label{I3}d_is_j^m=&s_j^md_{i-m}\hspace{0.5cm}&\ \mbox{if\qquad} i>j+m\\[3pt]
\label{I4}d_is_j^m=&s_j^{m-1}\hspace{0.5cm}&\ \mbox{if\qquad} j\leq i\leq j+m.
\end{eqnarray}
\end{lemma}

\vspace{0.3cm} Let $f:X\to Y$ be a simplicial map. A collection of simplices
\vspace{-0.2cm}
\begin{eqnarray*}
x_i&\in& X_{n-1}, \hspace{0.4cm} i\in I,\\
y&\in& Y_n\,,
\end{eqnarray*}
where $I\subseteq [n]$ is any subset, is said to be {\em $f$-compatible} whenever the following
equalities hold:
\vspace{-0.4cm}\begin{eqnarray*}d_ix_j=&d_{j-1}x_i&\hspace{0.4cm}\mbox{for all\qquad } i,j\in I, i<j,\\
d_iy=&fx_i~~~~ &\hspace{0.4cm} \mbox{for all\qquad }i\in I.\\[2pt]
\end{eqnarray*}

The map $f$ is said to be a {\em Kan fibration}
whenever for every given collection of $f$-compatible simplices
\begin{eqnarray*}x_0,\ldots,x_{k-1},x_{k+1},\ldots,x_n&\in& X_{n-1},\\
y&\in &Y_n\,,
\end{eqnarray*}
there is a
simplex $x\in X_n$ such that $d_ix=x_i$ for all $i\neq k$ and $fx=y$.
The next lemma (cf. \cite[Lemma 7.4]{may}) will be very useful in our development. For $I$ any finite
set, $|I|$ denotes its number of elements.
\begin{lemma}\label{1.1} Let $f:X\to Y$ be a Kan fibration. Suppose that there are given a subset $I\subseteq
[n]$ such that $1\leq |I|\leq n$ and an $f$-compatible family of simplices\linebreak
\vspace{-0.2cm}
\begin{eqnarray*}
x_i&\in& X_{n-1}, \hspace{0.4cm} i\in I,\\
y&\in& Y_n\,.
\end{eqnarray*}
Then, there exists $x\in X_n$ such that $d_ix=x_i$ for all $i\in I$ and $fx=y$.
\end{lemma}
\begin{proof} Suppose $|I|=r$. If $r=n$, the statement is true since $f$ is a Kan fibration. Hence
the statement holds for $n=1$. We now proceed by induction: Assume $n>1$ and the result holds for
$n'<n$ and assume $r<n$ and the result holds for $r'>r$. Taking $k=\mbox{max}\{i~|~i\in [n],
i\notin I\}$, we wish to find a simplex $x_k\in X_{n-1}$, such that the collection of simplices
\begin{eqnarray*}
x_i&\in& X_{n-1}, \hspace{0.4cm} i\in I\cup \{k\},\\
y&\in& Y_n\,,
\end{eqnarray*}
be $f$-compatible, since then an application of the induction hypothesis on $r$ gives the claim.
To find such an $x_k$, let $I'\subseteq [n-1]$ be the subset $$I'=\{i~|~i\in I, i<k\}\cup
\{i-1~|~i\in I, i>k\},$$and let
\begin{eqnarray*}
x'_i&\in& X_{n-2}, \hspace{0.4cm} i\in I',\\
y'&\in& Y_{n-1}\,,
\end{eqnarray*}
be the family of simplices defined by $x'_i=d_{k-1}x_i$ for $i\in I$, $i<k$, $x'_{i-1}=d_kx_i$ for
$i\in I$, $i>k$, and $y'=d_ky$. This family is $f$-compatible, whence the induction hypothesis
on $n-1$ gives the
required $x_k$, that is, a $(n-1)$-simplex of $X$ satisfying $d_ix_k=x'_i$ for $i\in I'$ and
$fx_k=y'$.
\end{proof}

\vspace{0.3cm} The category of {\em bisimplicial sets} is the category of functors
$X:\Delta^{op}\times \Delta^{op}\to \mbox{Set}$. It is often convenient to see a bisimplicial set
$X$ as a (horizontal) simplicial object in the category of (vertical) simplicial sets. If
$\alpha:[p]\to[p']$ and $\beta:[q]\to[q']$ are any two maps in $\Delta$, then we will write
$\alpha^{*_h}:X_{p',q}\to X_{p,q}$ and $\beta^{*_v}:X_{p,q'}\to X_{p,q}$ for the images
$X(\alpha,id)$ and $X(id,\beta)$ respectively. In particular, the horizontal and vertical face and
degeneracy maps are $d_i^h=(\delta_i)^{*_h}$, $d_i^v=(\delta_i)^{*_v}$, $s_i^h=(\sigma_i)^{*_h}$
and $s_i^v=(\sigma_i)^{*_v}$.

By composing with the diagonal functor $\Delta\to\Delta\times\Delta$, $[n]\mapsto ([n],[n])$, we
get the {\em diagonal functor} from bisimplicial sets to simplicial sets, which provides a
simplicial set $\diag X:[n]\mapsto X_{n,n}$, associated to each bisimplicial set $X$, whose face
and degeneracy operators are given in terms of those of $X$ by the formulas $d_i=d_i^hd_i^v$ and
$s_i=s_i^hs_i^v$, respectively.

\section{Proof of Theorem 1}
Let $p\geq 0$ be any fixed integer. Then, in order to  prove that the simplicial set map
$f_{p,*}:X_{p,*}\to Y_{p,*}$ is a Kan fibration,  suppose that, for some integers $q\geq 1$ and
$0\leq \ell \leq q$, there is given a collection of bisimplices
\begin{eqnarray*}x_0,\ldots,x_{\ell-1},x_{\ell+1},\ldots,x_q&\in& X_{p,q-1},\\
y&\in &Y_{p,q}\,,
\end{eqnarray*}
which is $f_{p,*}$-compatible, that is, such that $d_i^v x_j=d_{j-1}^v x_i$ for $0\leq i<j\leq q$,
$i\neq \ell\neq j$, and $d_i^v y=fx_i$ for $0\leq i\leq q$, $i\neq \ell$. We must therefore
find a bisimplex $ x\in X_{p,q}$ such that
\begin{equation}\label{r}\begin{array}{l}d_i^vx=x_i\hspace{0.6cm} \mbox{ for\quad }\  0\leq i\leq q, \ i\neq \ell,\\[4pt]fx=y
.\end{array}\end{equation}

\vspace{0.2cm}To do that, we start by considering the subset
$$ I\subseteq [p+q]$$
defined by
$$
I=\{i~|~0\leq i<\ell\}\cup\{p+i~|~\ell<i\leq q\},
$$
and the family of diagonal bisimplices
\begin{eqnarray}\label{3.2.1}
\nonumber\bar{x}_i&\in& (\diag X)_{p+q-1}, \hspace{0.6cm} i\in I,\\[-5pt]
\\[-5pt]
\nonumber\bar{y}&\in& (\diag Y)_{p+q}\,,
\end{eqnarray}
where
\begin{eqnarray}\label{-1}\bar{x}_i&=& (s_0^h)^{\ell-1}(s_p^h)^{q-\ell}
(s_{\ell-1}^v)^px_i \hspace{1cm}\mbox{for\quad }\hspace{0.3cm}0\leq i<\ell,
\\[4pt]\label{-2}
\bar{x}_{p+i}&=& (s_0^h)^{\ell}(s_p^h)^{q-\ell-1} (s_\ell^v)^p   x_i\hspace{1cm}\mbox{for\quad
}\hspace{0.3cm}\ell <i\leq q,
\\[4pt]\label{-3}
\bar{y}&=&(s_0^h)^{\ell}(s_p^h)^{q-\ell} (s_\ell^v)^py\,.
\end{eqnarray}

\vspace{0.2cm} The next verifications show that this family (\ref{3.2.1}) of diagonal bisimplices is
actually $\diag f$-compatible:

\vspace{0.3cm}\noindent- for $0\leq i< j<\ell$,
$$
\begin{array}{lcp{8cm}}
d_i\bar{x}_j&=&$d_i^hd_i^v(s_0^h)^{\ell-1}(s_p^h)^{q-\ell}(s_{\ell-1}^v)^{p} x_j$\\[7pt]
&=&$d_i^h(s_0^h)^{\ell-1}(s_p^h)^{q-\ell}d_i^v(s_{\ell-1}^v)^{p} x_j$\\[5pt]
&=&$(s_{0}^h)^{\ell-2}(s_{p}^h)^{q-\ell}(s_{\ell-2}^v)^{p}d_i^v x_j$\hfill{by (\ref{I4})}\\[7pt]
&=&$(s_{0}^h)^{\ell-2}(s_{p}^h)^{q-\ell}(s_{\ell-2}^v)^{p}d_{j-1}^v
x_i$\\[7pt]
&=&$(s_{0}^h)^{\ell-2}(s_{p}^h)^{q-\ell}d_{j-1}^v(s_{\ell-1}^v)^p x_i$\\[5pt]
&=&$d_{j-1}^h (s_0^h)^{\ell-1}(s_p^h)^{q-\ell}d_{j-1}^v(s_{\ell-1}^v)^p x_i$\hfill{by (\ref{I4})}\\[9pt]
&=&$d_{j-1}^h d_{j-1}^v(s_0^h)^{\ell-1}(s_p^h)^{q-\ell}(s_{\ell-1}^v)^p x_i$\\[7pt]&=&$d_{j-1}\bar{x}_i.$
\end{array}
$$

\vspace{0.3cm}\noindent- for $0\leq i<\ell< j\leq q$,
$$
\begin{array}{lcp{8cm}}
d_i\bar{x}_{p+j}&=&$d_i^hd_i^v(s_0^h)^{\ell}(s_p^h)^{q-\ell-1}(s_{\ell}^v)^{p} x_j$\\[7pt]
&=&$d_i^h(s_0^h)^{\ell}(s_p^h)^{q-\ell-1}d_i^v(s_{\ell}^v)^{p} x_j$\\[5pt]
&=&$(s_{0}^h)^{\ell-1}(s_{p}^h)^{q-\ell-1}(s_{\ell-1}^v)^{p}d_i^v x_j$\hfill{by (\ref{I4})}\\[7pt]
&=&$(s_{0}^h)^{\ell-1}(s_{p}^h)^{q-\ell-1}(s_{\ell-1}^v)^{p}d_{j-1}^v
x_i$\\[5pt]
&=&$(s_{0}^h)^{\ell-1}d_{p+j-\ell}^h(s_{p}^h)^{q-\ell}(s_{\ell-1}^v)^p
d_{j-1}^v
x_i$\hfill{by (\ref{I4})}\\[5pt]
&=&$d_{p+j-1}^h (s_0^h)^{\ell-1}(s_p^h)^{q-\ell}d_{p+j-1}^v(s_{\ell-1}^v)^p x_i$\hfill{by (\ref{I3})}\\[9pt]
&=&$d_{p+j-1}^h d_{p+j-1}^v(s_0^h)^{\ell-1}(s_p^h)^{q-\ell}(s_{\ell-1}^v)^p
x_i$\\[7pt]&=&$d_{p+j-1}\bar{x}_{i}.$
\end{array}
$$

\vspace{0.3cm}\noindent- for $\ell<i< j\leq q$,
$$
\begin{array}{lcp{8cm}}
d_{p+i}\bar{x}_{p+j}&=&$d_{p+i}^h(s_0^h)^{\ell}(s_p^h)^{q-\ell-1}d_{p+i}^v(s_{\ell}^v)^{p} x_j$\\[5pt]
&=&$(s_{0}^h)^{\ell}d_{p+i-\ell}^h(s_{p}^h)^{q-\ell-1}(s_{\ell}^v)^{p}d_i^v
x_j$\hfill{by (\ref{I3})}\\
[5pt]& =&$(s_{0}^h)^{\ell}(s_{p}^h)^{q-\ell-2}(s_{\ell}^v)^{p}d_{j-1}^v
x_i$\hfill{by (\ref{I4})}\\[5pt]
&=&$(s_{0}^h)^{\ell}d_{p+j-\ell-1}^h(s_{p}^h)^{q-\ell-1}(s_{\ell}^v)^p
d_{j-1}^v
x_i$\hfill{by (\ref{I4})}\\[5pt]
&=&$d_{p+j-1}^h (s_0^h)^{\ell}(s_p^h)^{q-\ell-1}d_{p+j-1}^v(s_{\ell}^v)^p x_i$\hfill{by (\ref{I3})}\\[9pt]
&=&$d_{p+j-1}^h d_{p+j-1}^v(s_0^h)^{\ell}(s_p^h)^{q-\ell-1}(s_{\ell}^v)^p
x_i$\\[7pt]&=&$d_{p+j-1}\bar{x}_{p+i}.$
\end{array}
$$

\vspace{0.3cm}\noindent- for $0\leq i<\ell$,

$$
\begin{array}{lcp{7cm}}
d_i\bar{y}&=&$d_i^hd_i^v(s_0^h)^{\ell}(s_p^h)^{q-\ell}(s_{\ell}^v)^{p} y$\\[7pt]
&=&$d_i^h(s_0^h)^{\ell}(s_p^h)^{q-\ell}d_i^v(s_{\ell}^v)^{p} y$\\[5pt]
&=&$(s_{0}^h)^{\ell-1}(s_{p}^h)^{q-\ell}(s_{\ell-1}^v)^{p}d_i^v y$\hfill{by (\ref{I4})}\\
[7pt]& =&$(s_{0}^h)^{\ell-1}(s_{p}^h)^{q-\ell}(s_{\ell-1}^v)^{p}fx_i $\\[7pt]
&=&$f(s_{0}^h)^{\ell-1}(s_{p}^h)^{q-\ell}(s_{\ell-1}^v)^{p}x_i$\\[7pt]
&=&$f\bar{x}_i.$
\end{array}
$$

\noindent and, finally,

\vspace{0.3cm}\noindent- for $\ell <i\leq q$,

$$
\begin{array}{lcp{7cm}}
d_{p+i}\bar{y}&=&$d_{p+i}^hd_{p+i}^v(s_0^h)^{\ell}(s_p^h)^{q-\ell}(s_{\ell}^v)^{p} y$\\[7pt]&=&
$d_{p+i}^h(s_0^h)^{\ell}(s_p^h)^{q-\ell}d_{p+i}^v(s_{\ell}^v)^{p} y$\\[5pt]
&=&$(s_{0}^h)^{\ell}d_{p+i-\ell}^h(s_{p}^h)^{q-\ell}(s_{\ell}^v)^{p}d_i^v
y$\hfill{by (\ref{I3})}\\[5pt]
&=&$(s_{0}^h)^{\ell}(s_{p}^h)^{q-\ell-1}(s_{\ell}^v)^{p}d_i^v y$\hfill{by (\ref{I4})}\\[5pt]
&=&$(s_{0}^h)^{\ell}(s_{p}^h)^{q-\ell-1}(s_{\ell}^v)^p fx_i$\\[7pt]&=&
$f(s_{0}^h)^{\ell}(s_{p}^h)^{q-\ell-1}(s_{\ell}^v)^p x_i$\\[5pt]
&=&$f\bar{x}_{p+i}.$
\end{array}
$$

\vspace{0.3cm}Then, since by hypothesis $\diag f:\diag X\to\diag Y$ is a Kan fibration, from Lemma
\ref{1.1} we get a diagonal bisimplex$$\bar{x}\in (\diag X)_{p+q}$$ such that
$d_i^hd_i^v\bar{x}=\bar{x}_i$ for $0\leq i<\ell$,  $d_{p+i}^hd_{p+i}^v\bar{x}=\bar{x}_{p+i}$ for
$\ell<i\leq q$ and $f\bar{x}=\bar{y}$.

\vspace{0.2cm}Now, using the bisimplex $\bar{x}$ we construct the bisimplex

$$x=(d_{p+1}^h)^{q-\ell}(d_0^h)^\ell(d_\ell^v)^p\,\bar{x}\in X_{p,q},$$

\noindent which, we claim, satisfies Relations (\ref{r}). Actually:

\vspace{0.3cm}\noindent- for $0<i<\ell$,

$$\begin{array}{lcp{9.4cm}}
d_i^v\,x&=&$d_i^v (d_{p+1}^h)^{q-\ell}(d_0^h)^\ell(d_\ell^v)^p\,\bar{x}$\\[5pt]
&=&$(d_{p+1}^h)^{q-\ell}(d_0^h)^\ell d_i^v(d_\ell^v)^p\,\bar{x}$\\[5pt]
&=&$(d_{p+1}^h)^{q-\ell}(d_0^h)^\ell(d_{\ell-1}^v)^pd_i^v\,\bar{x}$\\[5pt]
&=&$(d_{p+1}^h)^{q-\ell}(d_0^h)^{\ell-1}d_i^h(d_{\ell-1}^v)^pd_i^v\,\bar{x}$\hfill{by (\ref{I2})}\\[7pt]
&=&$(d_{p+1}^h)^{q-\ell}(d_0^h)^{\ell-1}(d_{\ell-1}^v)^pd_i^hd_i^v\,\bar{x}$\\[5pt]
&=&$(d_{p+1}^h)^{q-\ell}(d_0^h)^{\ell-1}(d_{\ell-1}^v)^p\bar{x}_i$\\[5pt]
&=&$(d_{p+1}^h)^{q-\ell}(d_0^h)^{\ell-1}(d_{\ell-1}^v)^p
(s_0^h)^{\ell-1}(s_p^h)^{q-\ell} (s_{\ell-1}^v)^px_i$\hfill{by (\ref{-1})}\\[7pt]
&=&$(d_{p+1}^h)^{q-\ell}(d_0^h)^{\ell-1} (s_0^h)^{\ell-1}(s_p^h)^{q-\ell}
(d_{\ell-1}^v)^p(s_{\ell-1}^v)^px_i$\\[7pt]&=&$x_i.$
\end{array}$$

\vspace{0.3cm}\noindent- for $\ell<i\leq q$,

$$\begin{array}{lcp{9.4cm}}
d_i^v\,x&=&$d_i^v (d_{p+1}^h)^{q-\ell}(d_0^h)^\ell(d_\ell^v)^p\,\bar{x}$\\[5pt]
&=&$(d_{p+1}^h)^{q-\ell}
(d_0^h)^\ell d_i^v(d_\ell^v)^p\,\bar{x}$\\[5pt]
&=&$(d_{p+1}^h)^{q-\ell}(d_0^h)^\ell(d_{\ell}^v)^pd_{p+i}^v\,\bar{x}$\hfill{by (\ref{I1})}\\[5pt]
&=&$(d_{p+1}^h)^{q-\ell-1}d_{p+i-\ell}^h(d_0^h)^{\ell}
(d_{\ell}^v)^pd_{p+i}^v\,\bar{x}$\hfill{by (\ref{I2})}\\[5pt]
&=&$(d_{p+1}^h)^{q-\ell-1}(d_0^h)^{\ell}d_{p+i}^h
(d_{\ell}^v)^pd_{p+i}^v\,\bar{x}$\hfill{by (\ref{I1})}\\[7pt]
&=&$(d_{p+1}^h)^{q-\ell-1}(d_0^h)^{\ell}
(d_{\ell}^v)^pd_{p+i}^hd_{p+i}^v\,\bar{x}$\\[7pt]
&=&$(d_{p+1}^h)^{q-\ell-1}(d_0^h)^{\ell}
(d_{\ell}^v)^p(s_0^h)^{\ell}(s_p^h)^{q-\ell-1} (s_\ell^v)^p   x_i$\hfill{by (\ref{-2})}\\[7pt]
&=&$(d_{p+1}^h)^{q-\ell-1}(d_0^h)^{\ell}
(s_0^h)^{\ell}(s_p^h)^{q-\ell-1}(d_{\ell}^v)^p (s_\ell^v)^p   x_i$\\[7pt]&=&$x_i.$
\end{array}$$

\noindent and, finally,

$$\begin{array}{lcp{9.4cm}}
fx&=&$ f(d_{p+1}^h)^{q-\ell}(d_0^h)^\ell(d_\ell^v)^p\,\bar{x}$\\
&=&$(d_{p+1}^h)^{q-\ell}(d_0^h)^\ell(d_\ell^v)^p f\,\bar{x}$\\[7pt]
&=& $(d_{p+1}^h)^{q-\ell}(d_0^h)^\ell(d_\ell^v)^p \bar{y}$\\[5pt]
&=&$
(d_{p+1}^h)^{q-\ell}(d_0^h)^\ell(d_\ell^v)^p (s_0^h)^{\ell}(s_p^h)^{q-\ell} (s_\ell^v)^py$\hfill{by (\ref{-3})}\\[7pt]
&=&$(d_{p+1}^h)^{q-\ell}(d_0^h)^\ell (s_0^h)^{\ell}(s_p^h)^{q-\ell}(d_\ell^v)^p
(s_\ell^v)^py$\\[7pt]&=&$ y.$
\end{array}$$

    Hence Theorem 1 is proved.

\vspace{0.2cm}    Let us stress that Theorem 1 not only says that, when $\diag f:\diag X\to \diag Y$ is a Kan
fibration,  all simplicial maps $f_{p,*}:X_{p,*}\to Y_{p,*}$ are Kan fibrations, but  that the
simplicial maps $f_{*,p}:X_{*,p}\to Y_{*,p}$ are Kan fibrations as well. This is a consequence of
the symmetry of the hypothesis: ``$\diag f$ is a Kan fibration", that is, the fact follows by
exchanging the vertical and horizontal directions.

\section{The converse of Theorem 1 is false}

 It is possible that, for a bisimplicial map $f:X\to Y$, all simplicial maps
$f_{p,*}$ and $f_{*,p}$, $p\geq 0$, be Kan fibrations and, however, $\diag\! f$ be not a Kan
fibration: Take $X$ to be the double nerve of an suitable double groupoid (i.e. a groupoid object
in the category of groupoids); then, all simplicial sets $X_{p,*}$ and $X_{*,p}$ are Kan complexes
(since they are nerves of groupoids) but, as we show below, the diagonal simplicial set $\diag X$
is not necessarily a Kan complex.

\vspace{0.2cm} We briefly recall some standard terminology about double groupoids; see for example
\cite{B-S} or \cite{K-S}. A double groupoid consists of objects $a$, $b$, ..., horizontal and vertical
morphisms between them
$$
a \overset{f}\to b,\hspace{0.5cm}\begin{array}{l}b\\\uparrow {\!_{g}}\\ a\end{array},\text{ ...},
$$
and squares $\sigma$, $\tau$, etc., of the form
$$
\xymatrix@C=0pt@R=0pt{c&&\ar[ll]_{f}b\\&\sigma&\\d\ar[uu]^{g'}&&\ar[ll]^{f'}a\ar[uu]_{g}}\,.
$$
These satisfy axioms such that the objects together with the horizontal morphisms as well as the
objects together with the vertical morphisms form groupoids. Furthermore, the squares form a
groupoid under both horizontal and vertical juxtaposition, and in the situation
$$
\xymatrix@C=0pt@R=0pt{\cdot&&\ar[ll]\cdot&&\cdot\ar[ll]\\
&\sigma& &\tau&\\ \ar[uu]\cdot&&\ar[ll]\cdot\ar[uu]&&\ar[uu]\cdot\ar[ll] \\
&\gamma& &\delta& \\ \ar[uu]\cdot&&\ar[ll]\cdot\ar[uu]&&\ar[ll]\cdot\ar[uu]}
$$
horizontal and vertical composition commute in the sense that $$(\sigma \cdot_h \tau)\cdot_v
(\gamma\cdot_h\delta) =(\sigma \cdot_v \gamma)\cdot_h (\tau\cdot_v\delta).$$ In addition to the
vertical identity morphisms, there are horizontal identity squares for each vertical morphism. Its
horizontal edges are the identity arrows of the horizontal groupoid of arrows. Similarly, there
are vertical identity squares for each horizontal morphism with vertical identity arrows. These
identity squares are compatible in the sense that vertical and horizontal identity squares for the
vertical and horizontal identity morphisms are the same.

Given a double groupoid $\mathbb{G}$, one can construct its bisimplicial nerve (or double nerve)
$\mbox{NN}\mathbb{G}$. A typical $(p,q)$ simplex of $\mbox{NN}\mathbb{G}$ is a subdivision of a
square of $\mathbb{G}$ as matrix of  $p\times q$ horizontally and vertically composable squares of
the form
$$
\xymatrix@C=0pt@R=0pt{\cdot&&\ar[ll]\cdot&&\cdot\ar[ll]&\cdots&\cdot&&\ar[ll]\cdot\\
&\sigma_{11}& &\sigma_{21}&&&&\sigma_{p1}&\\
\ar[uu]\cdot&&\ar[uu]\ar[ll]\cdot&&\ar[uu]\cdot\ar[ll]&\cdots&\ar[uu]\cdot&&\ar[uu]\ar[ll]\cdot  \\
&\sigma_{12}& &\sigma_{22}&&&&\sigma_{p2}&\\
\ar[uu]\cdot&&\ar[uu]\ar[ll]\cdot&&\ar[uu]\cdot\ar[ll]&\cdots&\ar[uu]\cdot&&\ar[uu]\ar[ll]\cdot\\
\vdots&&\vdots&&\vdots&&\vdots&&\vdots\\
\cdot&&\ar[ll]\cdot&&\cdot\ar[ll]&\cdots&\cdot&&\ar[ll]\cdot\\
&\sigma_{1q}& &\sigma_{2q}&&&&\sigma_{pq}&
\\\ar[uu]\cdot&&\ar[uu]\ar[ll]\cdot&&\ar[uu]\cdot\ar[ll]&\cdots&\ar[uu]\cdot&&\ar[uu]\ar[ll]\cdot
}
$$

The bisimplicial face operators are induced by horizontal and vertical composition of squares, and
degeneracy operators by appropriated identity squares. We picture $\mbox{NN}\mathbb{G}$ so that
the set of $(p,q)$-simplices is the set in the $p$-th row and $q$-th column. Thus, the $p$-th
column of $\mbox{NN}\mathbb{G}$, $\mbox{NN}\mathbb{G}_{p,*}$, is the nerve of the ``vertical"
groupoid whose objects are  strings of $p$ composable horizontal morphisms
$(a_0\overset{\hspace{5pt}f_1}\leftarrow a_1\overset{\hspace{5pt}f_2}\leftarrow\cdots
\overset{\hspace{5pt}f_p}\leftarrow a_p)$ and whose arrows are depicted as
$$\xymatrix@C=0pt@R=0pt{\cdot&&\ar[ll]\cdot&&\cdot\ar[ll]&\cdots&\cdot&&\ar[ll]\cdot\\
&\sigma_{1}& &\sigma_{2}&&&&\sigma_{p}&\\
\ar[uu]\cdot&&\ar[uu]\ar[ll]\cdot&&\ar[uu]\cdot\ar[ll]&\cdots&\ar[uu]\cdot&&\ar[uu]\ar[ll]\cdot
}$$ And, similarly, the $q$-th row, $\mbox{NN}\mathbb{G}_{*,q}$, is the ``horizontal" groupoid
whose objects are the length $q$ sequences of composable vertical morphisms in $\mathcal{C}$ and whose
arrows are sequences of $q$ vertically composable squares.  In particular,
$\mbox{NN}\mathbb{G}_{0,*}$ and $\mbox{NN}\mathbb{G}_{*,0}$ are, respectively, the nerves of the
groupoids of vertical and horizontal morphisms of $\mathbb{G}$.

\vspace{0.2cm} \noindent{\bf Example.} Let $A,B$ be subgroups of a finite group $G$, such that
$AB\neq BA$, where
$$AB=\{ab~|~a\in A,\,b\in B  \} \text{ and } BA=\{ba~|~a\in A,\,b\in B\}$$ ({\em for instance, we can
take $G=S_3$, $A=\{id, (1,2)\}$ and $B=\{id, (1,3)\}$}).

\vspace{0.2cm} Then, let $\mathcal{C}=\mathcal{C}(A,B)$ denote the double groupoid with only one
object, say ``$\cdot$",  with horizontal morphisms $\cdot\overset{a}\leftarrow \cdot$ the elements
$a\in A$, with vertical morphisms $\begin{array}{l}\cdot \\\uparrow {\!_{b}}\\ \cdot\end{array}$
are the elements $b$ of $B$, both with composition given by multiplication in the group, and whose
squares
$$
\xymatrix@C=3pt@R=3pt{\cdot &&\ar[ll]_{a} \cdot \\& & \\\cdot
\ar[uu]^{b'}&&\ar[ll]^{a'}\cdot\ar[uu]_{b}}
$$
are lists $(a,b,a',b')$ of elements $a,a'\in A$ and $b,b'\in B$ such that $ab=b'a'$. Compositions
of squares are defined, in the natural way, by
$$
(a,b,a',b')\cdot_h(a_1,b_1,a'_1,b)=(aa_1,b_1,a'a'_1,b'),
$$
$$
(a,b,a',b')\cdot_v(a',b_1,a'_1,b'_1)=(a,bb_1,a'_1,b'b'_1),
$$
and identities by
$$id^v(a)=(a,e,a,e),\ id^h(b)=(e,b,e,b),$$
where $e$ is the neutral element of the group.

\vspace{0.2cm} So defined, we claim that the associated simplicial set $\diag\mbox{NN}\mathcal{C}$
is not a Kan complex. To prove that, note that the hypothesis $AB\neq BA$ implies the existence of
elements $a\in A$ and $b\in B$ such that $ab$ cannot be expressed in the form $b'a'$, for any
$a'\in A$ and $b'\in B$. Now, the identity squares $\iota_b=id^h(b)$ and $\iota_{a}=id^v(a)$ can
be placed jointly as in the picture:
$$
\xymatrix@C=0pt@R=0pt{\cdot&&\ar[ll]_a\cdot&&\\
&\iota_{\!a}& &&\\ \ar@{=}[uu]\cdot&&\ar[ll]^{{a}}\cdot\ar@{=}[uu]&&\cdot\ar@{=}[ll] \\
&& &\iota_b& \\ &&\cdot\ar[uu]^b&&\ar@{=}[ll]\cdot\ar[uu]_b}
$$
This means that, regarded $\iota_a$ and $\iota_{b}$ as $1$-simplices of diagonal simplicial set
$\diag\mbox{NN}\mathcal{C}$, they are compatible, in the sense that
$$d^h_0d^v_0\iota_{a}=\cdot=d^h_1d^v_1\iota_b.$$ Therefore, if we assume that $\diag\mbox{NN}\mathcal{C}$
is a Kan complex, we should find a diagonal bisimplex, say $x\in \mbox{NN}\mathcal{C}_{2,2}$, such
that $d^h_0d^v_0x=\iota_b$ and $d^h_2d^v_2x=\iota_{a}$. This $(2,2)$-simplex $x$ of the nerve of
the double groupoid, should be therefore of the form
$$
\xymatrix@C=0pt@R=0pt{\cdot &&\ar[ll]_a\cdot&&\cdot\ar[ll]\\
&\iota_{a}& &&\\ \ar@{=}[uu]\cdot&&\ar[ll]^{\ {a}}\cdot\ar@{=}[uu]&&\cdot\ar@{=}[ll]\ar[uu] \\
&& &~\iota_b& \\\ar[uu]^{a'}\cdot &&\ar[ll]^{b'}\cdot\ar[uu]_b&&\ar@{=}[ll]a\ar[uu]_b}
$$
what is clearly impossible. Hence, the diagonal of the bisimplicial nerve of $\mathcal(A,B)$ is
not a Kan complex, in spite of the fact that this bisimplicial set is a pointwise Kan complex in
both vertical and horizontal directions.

\vspace{0.2cm}

For a last comment,  we shall point out a fact suggested by a referee:  the
property proved in Theorem 1 for diagonal fibrations is not true for diagonal trivial fibrations.
For instance: Let $EG$ be the universal cover of a non-trivial discrete group $G$, that is, the simplicial set with $EG_n=G^{n+1}$ and faces given by $d_i(x_0,...,x_n)=(x_0,...,x_{i-1},x_{i+1},...,x_n)$, $0\leq i\leq n$,  and let $X=EG\otimes
EG$, the bisimplicial set defined by $X_{p,q}=EG_p\times EG_q$. Then, $\diag X$ is
contractible to a one point Kan complex and, however, each simplicial set $X_{*,q}=EG\times G^q$ is
homotopy equivalent to the constant simplicial set $G^q$, and therefore it is not contractible. In
conclusion, $X\to pt$ is a diagonal trivial fibration but it is not a pointwise trivial fibration.

$\Delta$.
\end{document}